\documentclass[11pt]{amsart}
\usepackage{amsmath}
\usepackage{amsthm}
\usepackage{amssymb}
\usepackage{graphicx}

\newtheorem{theorem}{Theorem}[section]
\newtheorem{lemma}[theorem]{Lemma}
\newtheorem{prop}[theorem]{Proposition}

\newtheorem{defn}[theorem]{Definition}
\newtheorem{ex}[theorem]{Example}

\newcommand{\CC}{\mathbb{C}}

\newcommand{\PP}{\mathbb{P}}

\begin{document}
\title{Hilbert polynomials and module generating degrees}

\author{Roger Dellaca}
\address{University of California Irvine, Mathematics Department, Irvine, CA 92697, USA}
\email{rdellaca@uci.edu}

\begin{abstract}
We establish a form of the Gotzmann representation of the Hilbert polynomial based on rank and generating degrees of a module, which allow for a generalization of  Gotzmann's Regularity Theorem. Under an additional assumption on the generating degrees, the Gotzmann regularity bound becomes sharp. An analoguous modification of the Macaulay representation is used along the way, which generalizes the theorems of Macaulay and Green, and Gotzmann's Persistence Theorem. 
\end{abstract} 


\maketitle


\section{Introduction}

Several results on Hilbert functions and Hilbert polynomials can be interpreted in terms of extremal behavior achieved by lex ideals. In particular, Macaulay's Theorem and Green's Hyperplane Restriction Theorem give inequalities on the Hilbert function which are sharp in the case of a quotient by a lex ideal, in all degrees larger than the maximal generating degree of the ideal. Also, Gotzmann's Regularity Theorem gives a bound on the Castelnuovo-Mumford regularity of a saturated ideal, which is achieved by the saturated lex ideal. Gasharov extended the theorems of Macaulay and Green to modules \cite{Gas97}, and this author extended Gotzmann's Regularity Theorem to globally generated coherent sheaves \cite{Del15}; however, these results are not sharp as seen in \cite[Remark 4.1]{Gas97} and Example \ref{ex-not-sharp} below.

Hulett gave a generalization of Macaulay's Theorem which is a sharp inequality in case of lex modules; see \cite{Hul95} and \cite{BE00}. Additionally, Greco generalized Green's theorem in such a way that the result is sharp for lex modules \cite{Gre15}. Both of these use the same modification of the Macaulay representation and similar modifications of the Macaulay and Green transformations. In this article, we will see a modification of the  Gotzmann representation which accounts for the rank and generating degrees of a module. Based on this representation, a generalization of Gotzmann's Regularity Theorem is given which is sharp under an additional assumption on the generating degrees (see Theorem \ref{thm-rankadj-Gotzmann-reg}).

In order to prove the generalized theorem on regularity, we will use a representation of the Hilbert function that is analogous to the rank-and-degree adjusted Gotzmann representation of the Hilbert polynomial. This is less granular than Hulett's representation, but it uses the same number of repeated terms in all degrees. We will prove analogous generalizations of the theorems of Macaulay and Green. In addition, a generalized Gotzmann's Persistence Theorem will be proven; this allows us to embed the Quot scheme $\mathrm{Quot}^P_{\mathcal{O}_{\mathbb{P}^n}}$ into a Grassmannian using the construction given in \cite{Del15}; we will see that the refinement here allows embedding into a substantially smaller dimensional Grassmannian. 

The structure of the paper is as follows. Section 2 gives background definitions and results. Section 3 establishes Macaulay and Gotzmann representations adjusted for rank and generating degrees. Section 4 applies the rank-and-degree adjusted Gotzmann to show a sharp inequality on the first and second Chern classes of a globally generated coherent sheaf, improving the results given in \cite{Del15}. Section 5 proves the generalizations of Macaulay, Green and Gotzmann Regularity Theorems, and shows that the Gotzmann regularity bound is achieved by the  saturated lex submodule with the same Hilbert polynomial, under assumptions on the generating degrees; these are satisfied, for example, by quotients of a free module generated in degree zero. Section 6 proves the generalization of Gotzmann's Persistence Theorem and gives some examples to show the improvement in the regularity bound.


\section{Background}

Let $k$ be a field and $S=k[x_0, \ldots, x_n]$. Assume $F$ is a graded free $S$-module $F = Se_1 + \cdots + Se_m$ with $\deg(e_i) = f_i$ and $f_1 \leq \cdots \leq f_m$, and that $N$ is a graded submodule of $F$, with $M=F/N$ having rank $r$.

The \emph{Hilbert function} of $M$ is $H(M,d) = \dim_k(M_d)$, the dimension of the degree-$d$ part of $M$ as a $k$-vector space. For $d \gg 0$, the Hilbert function becomes a polynomial $P_M(d)$, the \emph{Hilbert polynomial} of $M$.

\begin{defn}Given $a, d \in \mathbb{N}$, the \emph{$d$th Macaulay representation of $a$} is the unique expression
\[
a = \binom{k(d)}{d} + \binom{k(d-1)}{d-1} + \cdots + \binom{k(\delta)}{\delta},
\]
with $\delta \in \mathbb{Z}$, satisfying $k(d) > \cdots > k(\delta) \geq \delta > 0$. Given this representation, the \emph{$d$th Macaulay transformation of $a$} is
\[
a^{\langle d \rangle} = \binom{k(d) + 1}{d + 1} + \binom{k(d-1) + 1}{d} + \cdots + \binom{k(\delta) + 1}{\delta + 1},
\]
 and the \emph{Green transformation} is
\[
a_{\langle d \rangle} = \binom{k(d) - 1}{d} + \binom{k(d-1) - 1}{d-1} + \cdots + \binom{k(\delta) - 1}{\delta}.
\]
\end{defn}

\begin{theorem} \label{thm-Green-Macaulay-Gotz-Persistence} Let $R = S/I$ be a graded $k$-algebra and $d \geq 1$ and integer. 
\begin{enumerate}
\item (Green's Hyperplane Restriction Theorem) For a general linear form $h$,
\[
H(R/hR, d) \leq H(R,d)_{\langle d \rangle}
\] 
\item (Macaulay)\cite{Mac27} $H(R,d+1) \leq H(R,d)^{\langle d \rangle}$
\item (Gotzmann's Persistence Theorem)\cite{Got78} If $I$ is generated in degree at most $d$ and $H(S/I,d+1) = H(S/I,d)^{\langle d \rangle}$, then 
\[
H(R,s+1) = H(R,s)^{\langle s \rangle}
\]
for all $s \geq d$.
\end{enumerate}
\end{theorem}

\proof 
\begin{enumerate}
\item See \cite[Theorem 1]{Gre89};
\item See \cite[Theorem 4.2.10]{BH98};
\item See \cite{Gre89}. $\qed$
\end{enumerate}

\begin{defn}
Given a polynomial $P(d) \in \mathbb{Q}[d]$, a \emph{Gotzmann representation} of $P$ is a binomial expansion
\[
P(d) = \binom{d + a_1}{d} + \binom{d+a_2 - 1}{d-1} + \cdots + \binom{d + a_s - (s-1)}{d-(s-1)},
\]
with $a_1, \cdots, a_s \in \mathbb{Z}$ and $a_1 \geq \cdots \geq a_s \geq 0$. The number of terms $s$ in the representation is the \emph{Gotzmann number} of $P$. 
\end{defn}

\begin{defn}
A finitely-generated graded $S$-module is $d$ regular if 
\[
d \geq \max\{e | H_{\mathfrak{m}}^i(M)_e \neq 0\} + i \text{ for all } i\geq 0.
\]
The Castelnuovo-Mumford regularity of $M$ is the smallest such $d$.
\end{defn}

\begin{theorem}(Gotzmann's Regularity Theorem \cite{Got78})
The Hilbert polynomial of a standard graded algebra has a unique Gotzmann representation. Furthermore, if $R=S/I$ is a standard graded algebra with Gotzmann number $s$, then $\mathrm{sat}(I)$ is $s$-regular.
\end{theorem}

\proof See \cite[Corollary B.31]{Vas04} for the first statement and \cite[Theorem B.33]{Vas04} for the second statement. $\qed$

Theorem \ref{thm-Green-Macaulay-Gotz-Persistence} was generalized to modules by \cite{Gas97}:

\begin{theorem} \label{Gas97thm}
Assume $F$ is a free graded $S$-module with $l$ the maximal degree of its generators,
$N$ a submodule of $F$ and $M = F /N$. Let $x$ be a general element of $S_1$,
$S' = S/(x)$, and $M' = F /(N + xF)$. For each pair $(p, d)$ such that $p \geq 0$
and $d \geq p + l + 1$, we have:
\begin{enumerate}
\item $H(M, d + 1) \leq H(M, d)^{\langle d - l - p \rangle}$;
\item $H(M',d) \leq H(M,d)_{\langle d - l - p \rangle}$;
\item If $N$ is generated in degree at most $d$ and
$H(M, d + 1) = H(M, d)^{\langle d - l - p \rangle}$,
then
$H(M, d + 2) = H(M, d + 1)^{\langle d+1-l-p \rangle}$.
\end{enumerate}
\end{theorem}

\proof See \cite[Theorem 4.2]{Gas97}. $\qed$

Gotzmann's regularity theorem was generalized to globally generated coherent sheaves by \cite{Del15}:

\begin{theorem} \label{thm-Gotz-globgen} If $f_m \leq 0$ and $N$ is a submodule of $F$, then writing $M=F/N$, we have $P_{M}(d)$ has a unique Gotzmann representation. Furthermore, if $P_M(d)$ has Gotzmann number $s$, then $\mathrm{sat}(N)$ is $s$-regular.
\end{theorem}

\proof See \cite[Propositions 3.1 and 4.1]{Del15}. $\qed$

The bounds in Macaulay, Green, and Gotzmann regularity are sharp, achieved by lex ideals (see \cite[Proposition 4.2.8]{BH98} for Macaulay; \cite[Proposition 5.5.23]{KR05} for Green; and \cite[Proposition C.24]{Iar99} for Gotzmann regularity). However, these bounds are no longer sharp in Theorems \ref{Gas97thm} and \ref{thm-Gotz-globgen}, as shown by the next example.

\begin{ex} \label{ex-not-sharp}
Let $S=\CC[x,y]$ and $F=Se_1 + Se_2 + Se_3$, with $\deg(e_i)=0$ for each $i$. Take the saturated lex module $L=Se_1$. Then $H(F/L, 1) = 4 = \binom{4}{1}$, and $L$ has no generators in degree $2$, but $H(F/L,2) = 6 < 10 = 4^{\langle 1 \rangle}$. Next, for a linear form $h$, taking $F'/L'$ the image of $F/L$ in $S'=S/(h)$, we have $H(F'/L', 2) = 2 > 0 = 4_{\langle 1 \rangle}$. Finally, we note that $F/L$ is $0$-regular, but 
\[
P_{F/L}(d) = \binom{d+1}{1} + \binom{d}{1} + 1
\]
has Gotzmann number $3$.
\end{ex}

Hulett's representation of the Hilbert function of a module takes into consideration the generating degrees, and the Macaulay transformation is adjusted accordingly; see \cite[Theorem 3.2]{BE00} and \cite[Theorem 3.18]{Gre15} for the sharp generalizations of Macaulay and Green, respectively. In the next section, we give the rank-and-degree adjusted representations that will allow us to prove a sharp result on Gotzmann regularity for modules later.


\section{Macaulay and Gotzmann representations adjusted for rank and degree}

Assume $M=F/N$ as defined in the beginning of the last section. The following proposition establishes representations of the  Hilbert function and Hilbert polynomial that are adjusted by the rank and generating degrees.

\begin{prop} \label{prop-rankadj-rep}
The Hilbert function of $M$ can be written 
\[
H(M,d) = \sum_{i=m-r+1}^m\binom{d-f_i + n}{n} + \rho_d,
\]
with $0 \leq \rho_d \leq \sum_{i=1}^{m-r}\binom{d-f_i+n}{n}$ for each $d$. 
If $f_{m-r} \leq 0$, then the Hilbert polynomial of $M$ can be written
\[
P_M(d) = \sum_{i=m-r+1}^m\binom{d-f_i + n}{n} + Q_M(d),
\]
with $Q_M(d)$ having a unique Gotzmann representation. 

\end{prop}

\proof
By \cite{Big93},\cite{Hul95}, and \cite{Par96}, there exists a lex submodule $L \subset F$ such that $H(M,d) = H(F/L, d)$. Note that the rank of $F/L$ is the same as the rank of $M$. It follows that 
\[
F/L = S(-f_{1})/I_{1} \oplus \cdots \oplus S(-f_{m-r})/I_{m-r} \oplus \bigoplus_{i=m-r+1}^m{S(-f_i)} 
\]
for some monomial ideals $I_1, \ldots, I_r$ of $S$. The existence of the first  representation is immediate. 

For the second representation,  apply Theorem \ref{thm-Gotz-globgen} to 
\[
S(-f_{r+1})/I_{r+1} \oplus \cdots \oplus S(-f_m)/I_m
\]
to give the result. $\qed$

\begin{defn}
The Hilbert polynomial representation in the last proposition will be called a \emph{rank-and-degree adjusted Gotzmann representation} of $P_M$, and the the \emph{rank-and-degree adjusted Gotzmann number} of $P_M$ will be the Gotzmann number of $Q_M$.
\end{defn}

If $M=S/I$ is a quotient of $S$, then $r=1$ and $m=f_1=0$, and in this case we specialize to the Gotzmann representation for the Hilbert polynomial of $S/I$.

\begin{ex} \label{ex-Gotzmannrep}
Set $S=\CC[x,y,z]$ and $F=Se_1 + Se_2 + Se_3$, where the generating degrees of $F$ are $f_1=f_2=-1, f_3=0$. Assume $M=F/Se_3$. The Hilbert polynomial of $M$ is $P_M(d) = 2\binom{d+3}{2}$. The rank-and-degree adjusted Gotzmann representation is
\[
P_M(d) = \binom{d+3}{2} + \binom{d+2}{2} + \binom{d+1}{d} + 1,
\]
and the standard Gotzmann representation is
\[
P_M(d) = \binom{d+2}{d} + \binom{d+1}{d-1} + \binom{d-1}{d-2} + \binom{d-2}{d-3} + \binom{d-3}{d-4} + 11.
\]
The rank-and-degree adjusted Gotzmann number 2, and the standard Gotzmann number is 16.
\end{ex}


\section{Sharp bounds on first and second Chern classes}

As a first application, we will refine the inequality on Chern classes of a rank $r=2$ globally generated coherent sheaf on $\mathbb{P}^3$ given in \cite{Del15}. This bound will be sharp, and applies to all positive ranks $r$ and all dimensions $n$ of $\mathbb{P}^n$.

\begin{prop}
Let $E$ be a rank $r \geq 1$ globally generated coherent sheaf on $\mathbb{P}^n$ for an integer $n \geq 0$. Then the first and second Chern classes of $E$ satisfy the inequality
\[
c_2 \leq c_1^2.
\]
\end{prop}

\proof
If $n\leq 1$, then $c_2 = 0$, and $c_1 \geq 0$ for any such globally generated sheaf on $\mathbb{P}^n$, so assume that $n \geq 2$. 
We will only use the top 3 terms of the Hilbert polynomial, so it is sufficient to write the Todd class and Chern character as
\[
Td(\mathbb{P}^n) = 1 + \frac{n+1}{2}H + \frac{(n+1)(3n+2)}{24}H^2 + \cdots
\]
and
\[
ch(E(d)) = \sum_{i=0}^m{(-1)^i ch(F_i(d))},
\]
where
\[
0 \to F_m \to \cdots \to F_0 \to E
\]
is a resolution of $E$ by locally free sheaves, and 
\begin{align*}
ch(F_i(d)) = & \frac{(\alpha_{i,1} + d)^n + \cdots + (\alpha_{i,r_i} + d)^n}{n!}H^n \\
   & \quad + \frac{(\alpha_{i,1} + d)^{n-1} + \cdots + (\alpha_{i,r_i} + d)^{n-1}}{(n-1)!}H^{n-1} \\
   & \quad + \frac{(\alpha_{i,1} + d)^{n-2} + \cdots + (\alpha_{i,r_i} + d)^{n-2}}{(n-2)!}H^{n-2} + \cdots
\end{align*}

Now, we compute the top 3 terms of the Hilbert polynomial using Hirzebruch-Riemann-Roch:
\begin{align*}
P_E(d) & = r \frac{d^n}{n!} + \left(\frac{r(n+1)}{2} + c_1 \right)\frac{d^{n-1}}{(n-1)!} \\
  & \quad + \left(\frac{c_1^2-2c_2 + (n+1)c_1}{2} + \frac{r(n+1)(3n+2)}{24}\right)\frac{d^{n-2}}{(n-2)!} + \cdots,
\end{align*}
where $c_i = c_i(E)$.

Write the rank-and-degree adjusted Gotzmann representation of $P_E$ as
\[
P_E(d) = P_n + P_{n-1} + \cdots + P_0,
\]
where $P_i$ is the sum of binomial coefficients in the Gotzmann representation of degree $i$. Note that 
\begin{align*}
P_n & = r{d+n \choose n} \\
   & = r\frac{d^n}{n!} + \frac{r(n+1)}{2}\frac{d^{n-1}}{(n-1)!} + \frac{r(n+1)(3n+2)}{24}\frac{d^{n-2}}{(n-2)!} + \cdots,
\end{align*}
where the third term is determined by 
\[
\sum_{1\leq i<j\leq n}ij = \frac{(n-1)n(n+1)(3n+2)}{24}
\]
using a straightforward induction argument.

Then
\[
P_E - P_n = c_1 \frac{d^{n-1}}{(n-1)!} + \left(\frac{c_1^2 - 2c_2 + (n+1)c_1}{2}\right)\frac{d^{n-2}}{(n-2)!} + \cdots
\]
If $c_1 = 0$, then $P_E - P_n =  -c_2 \frac{d^{n-2}}{(n-2)!} + \cdots$, which requires $c_2 \leq 0$ for there to be a Gotzmann representation; otherwise, 
\begin{align*}
P_{n-1} & = {d+n-1 \choose n-1} + {d+n-2 \choose n-1} + \cdots + {d+n-c_1 \choose n-1} \\
   & = c_1\frac{d^{n-1}}{(n-1)!} + \left(\frac{(n+1)c_1 - c_1^2}{2}\right)\frac{d^{n-2}}{(n-2)!} + \cdots
\end{align*}

So we have
\[
P_E - P_n - P_{n-1} = (c_1^2 - c_2)\frac{d^{n-2}}{(n-2)!} + \cdots,
\]
with a non-negative leading coefficient for there to be a rank-adjusted Gotzmann representation, therefore $c_2 \leq c_1^2$. $\qed$

\begin{ex}
Let $S=\mathbb{C}[x_0, \ldots, x_n]$, with $n \geq 2$, let $a>0$ be an integer, and take 
\[
N=(x_0^a,x_1^a)S,
\]
and let $\mathcal{E}$ be the sheaf associated to $N$. Then $\mathcal{E}(a)$ is globally generated of rank 1, and from the resolution
\[
0 \to S(-a) \to S^2 \to N(a) \to 0
\]

it follows that $\mathcal{E}(a)$ has Chern polynomial
\[
c_t(\mathcal{E}(a)) = \frac{1}{1-at} = 1 + at + a^2t^2 + \cdots,
\]
with $c_1 = a$ and $c_2 = a^2$ achieving the strict bound in the theorem. For $r>1$, take a direct sum of $E$ with $\mathcal{O}^{r-1}$.
\end{ex}


\section{Rank and degree-adjusted Macaulay, Green, and Gotzmann Regularity}

We will extend Macaulay's theorem, Green's Theorem, and Gotzmann's Regularity Theorem in this section; Gotzmann's Persistence Theorem will be extended in the next section. As before, assume $F=Se_1 + \cdots + Se_m$ with $\deg(e_i) = f_i$ and $f_1 \leq \cdots \leq f_m$.

First we will require a lemma giving more properties of Green and Macaulay transformations, and a lemma on the behavior of regularity in exact sequences.

\begin{lemma} \label{lem-Gasharov-transforms}
Assume $a$, $b$ and $d$ are positive integers.
\begin{enumerate}
\item $a_{\langle d \rangle} + b_{\langle d \rangle} \leq (a+b)_{\langle d \rangle}$
\item $a^{\langle d \rangle} + b^{\langle d \rangle} \leq (a+b)^{\langle d \rangle}$
\item $a_{\langle d+1 \rangle} \leq a_{\langle d \rangle}$
\item $a^{\langle d+1 \rangle} \leq a^{\langle d \rangle}$
\end{enumerate}
\end{lemma}

\proof See \cite[Lemmas 4.4 and 4.5]{Gas97}. $\qed$

\begin{lemma} \label{lem-reg-exactseq}
Suppose 
\[
0 \to M' \to M \to M'' \to 0
\]
is an exact sequence of finitely-generated graded $S$-modules.
\begin{enumerate}
\item $\mathrm{reg}(M) \leq \mathrm{max}(\mathrm{reg}(M'), \mathrm{reg}(M''))$;
\item $\mathrm{reg}(M') \leq \mathrm{max}(\mathrm{reg}(M), \mathrm{reg}(M'')+1)$.
\end{enumerate}
\end{lemma}

\proof See \cite[Corollary 18.7]{Pee11}. $\qed$

We are ready to generalize Macaulay's and Green's Theorems.

\begin{prop} \label{prop-rankadj-Green-Macaulay}
Assume that $N$ is a submodule of $F$ such that $M=F/N$ has rank $r$, and Hilbert function
\[
H(M,d) = \sum_{i=m-r+1}^m\binom{d-f_i + n}{n} + \rho_d.
\]
\begin{enumerate}
\item For all $d \geq f_{m-r}+1$,
\[
H(M,d+1) \leq \sum_{i=m-r+1}^m\binom{d+1-f_i+n}{n} + \rho_d^{\langle d-f_{m-r} \rangle}.
\] 

\item For a general element $h \in S_1$, writing $F'=F/hF$, $M'=F/(N+hF)$ and $S'=S/hS$, we have
\[
H(M',d) \leq \sum_{i=m-r+1}^m\binom{d-f_i + n-1}{n-1} + (\rho_d)_{\langle d - f_{m-r} \rangle}
\]
for all $d \geq f_{m-r}+1$, where $H(M',d) = \dim_{S'}(M'_d)$.
\end{enumerate}
\end{prop}

\proof
The lex module with the same Hilbert function can be written 
\[
L = \bigoplus_{i=1}^{m-r}{I_i}e_i
\]
for some lex ideals $I_1, \ldots, I_r$ of $S$. Note that 
\begin{align*}
H(M,d) & = H(F/L,d) \\
    & = \sum_{i=m-r+1}^m\binom{d-f_i + n}{n} + H(S/I_1,d-f_1) + \cdots H(S/I_{m-r},d-f_{m-r}).
\end{align*}

Then for $d \geq f_{m-r}+1$, 
\begin{align*}
H(M,d+1) & = H(F/L,d+1) \\
    & \leq \sum_{i=m-r+1}^m\binom{d+1-f_i+n}{n} \\
		& \qquad + H(S/I_1,d-f_1)^{\langle d - f_1 \rangle} + \cdots H(S/I_{m-r},d-f_{m-r})^{\langle d-f_{m-r} \rangle} \\
		& \leq \sum_{i=m-r+1}^m\binom{d+1-f_i+n}{n} \\
		& \qquad + H(S/I_1,d-f_1)^{\langle d - f_{m-r} \rangle} + \cdots H(S/I_{m-r},d-f_{m-r})^{\langle d-f_{m-r} \rangle} \\
		& \leq \sum_{i=m-r+1}^m\binom{d+1-f_i+n}{n} + \rho_d^{\langle d-f_{m-r} \rangle}
\end{align*}
where the first inequality is by Macaulay's Theorem, the second and third inequalities are by Lemma \ref{lem-Gasharov-transforms}, thus proving the first statement.

For the second statement, write $I_i'$ and for the image of $I_i$ in $S'$ and $L'$ for the image of $L$ in $F'$. By \cite[Theorem 3.18]{Gre15}, we have $H(M',d) \leq H(F'/L', d)$, and $H(S'/I_i', d) = H(S/I, d)_{\langle d \rangle}$ for each lex ideal $I_i$ and for all $d > 0$ by \cite[Proposition 5.5.23]{KR05}. Then for $d \geq f_{m-r}+1$, with the last two inequalities by Lemma \ref{lem-Gasharov-transforms},
\begin{align*}
H(M',d) & \leq H(F'/L',d) \\
    & = \sum_{i=m-r+1}^m\binom{d-f_i+n-1}{n-1} \\
		& \qquad + H(S/I_1,d-f_1)_{\langle d - f_1 \rangle} + \cdots H(S/I_{m-r},d-f_{m-r})_{\langle d-f_{m-r} \rangle} \\
		& \leq \sum_{i=m-r+1}^m\binom{d-f_i+n-1}{n-1} \\
		& \qquad + H(S/I_1,d-f_1)_{\langle d - f_{m-r} \rangle} + \cdots H(S/I_{m-r},d-f_{m-r})_{\langle d-f_{m-r} \rangle} \\
		& \leq \sum_{i=m-r+1}^m\binom{d-f_i+n-1}{n-1} + (\rho_d)_{\langle d-f_{m-r} \rangle}. \qed
\end{align*}

Next, the rank and degree-adjusted form of Gotzmann regularity is proved.

\begin{theorem} \label{thm-rankadj-Gotzmann-reg}
Assume $F=Se_1 + \cdots + Se_m$ with $\deg(e_i) = f_i$ and $f_1 \leq \cdots \leq f_m$, such that $f_{m-r}\leq 0$, and $N$ is a submodule of $F$ such that $M=F/N$ has rank $r$, and Hilbert polynomial
\[
P_M(d) = \sum_{i=m-r+1}^m\binom{d-f_i + n}{n} + Q(d)
\]
where $Q_M(d)$ has Gotzmann number $s$. Then the saturation of $N$ is $\mathrm{max}(s,f_m)$-regular. Furthermore, if $f_{m-r}=0$ and $s \geq f_m$, then this bound is achieved by the lex submodule.
\end{theorem}

\proof
The proof is the same as the proof for modules in \cite{Del15}, by induction on $\dim(S) - 1$, and will be sketched.
From a saturated module $N$ with $M=F/N$ and a general linear form $h$, we obtain $M'=M/hM(-1)$ which still has rank $r$ since $h$ is $M$-regular, and satisfies the hypotheses of Proposition \ref{prop-rankadj-rep}, so the uniqueness of the Gotzmann representation for $Q_{M'}(d)$ gives
\[
P_{M'}(d) = \sum_{i=m-r+1}^m\binom{d-f_i + n-1}{n-1} + \binom{d+a_1 - 1}{a_1-1} + \cdots + \binom{d+a_t - (t-1)}{a_t-1}
\]
for some $t \leq s$, where 
\[
Q_M(d) = \binom{d+a_1}{a_1} + \cdots + \binom{d+a_s-(a_s-1)}{a_s}
\]
is the Gotzmann representation of $Q_M(d)$. The same concluding argument as in \cite{Del15}, using Proposition \ref{prop-rankadj-Green-Macaulay} in place of Theorem \ref{Gas97thm}, shows that $M$ is $(s-1)$-regular, and since $F$ is $f_m$-regular, by Lemma \ref{lem-reg-exactseq} (2), $N$ is $\mathrm{max}(s,f_m)$-regular.

To show sharpness when $f_{m-r}=0$ and $s \geq f_m$, consider a given Hilbert polynomial
\[
P(d) = \sum_{i=m-r+1}^m\binom{d-f_i + n}{n} + Q(d)
\]
with $Q(d)$ having Gotzmann number $s$. The saturated lex module $L$ with $F/L$ having the same Hilbert polynomial is
\[
L=S(-f_1) \oplus \cdots \oplus S(-f_{m-r-1}) \oplus L_{m-r}.
\]
for a lex ideal $L_{m-r}$ with $P_{S/L_{m-r}}(d) = Q(d)$.

Since $f_{m-r}=0$ and $\mathrm{reg}(L_{m-r})=s$, it follows that $\mathrm{reg}(L) = s$.
 $\qed$

In particular, rank-adjusted Gotzmann regularity is sharp for modules generated in degree zero. This is in contrast to Example \ref{ex-not-sharp} and \cite[Example 6.3]{Del15}, where no submodule achieves the regularity given by the standard Gotzmann number.


\section{Rank and degree-adjusted Gotzmann persistence}

It remains to prove rank-and-degree adjusted Gotzmann Persistence. The proof will be similar to \cite{GMP11}.

\begin{theorem}\label{thmPersistence} 

Assume (in addition to previous assumptions on $F$, $N$, and $M=F/N$) that $N$ is generated in degree at most $d$ and $d \geq f_{m-r}+1$. If 
\[
H(M,d+1) = \sum_{i=m-r+1}^m\binom{d+1-f_i+n}{n} + \rho_d^{\langle d-f_{m-r} \rangle},
\]
then
\[
H(M,d+2) = \sum_{i=m-r+1}^m\binom{d+2-f_i+n}{n} + \rho_{d+1}^{\langle d-f_{m-r}+1 \rangle}.
\]
\end{theorem}

The proof follows from extensions of some other results to modules, which will be given in two lemmas. We require a definition before the first lemma.

\begin{defn}
Suppose the module $M$ has graded Betti numbers $\beta_{i,j}$. A \emph{consecutive cancellation} is the process of choosing $i,j$ such that $\beta_{i,j}$ and $\beta_{i+1,j}$ are positive, and replacing $\beta_{i,j}, \beta_{i+1,j}$ with $\beta_{i,j}-1, \beta_{i+1,j}-1$ respectively.
\end{defn}

The following Lemma is the extension of \cite[Theorem 1.1]{Pee04} to modules. Recall that a \emph{weight ordering} $<_{\mathbf{w}}$ with weight vector $\mathbf{w}=(w_0, \ldots w_n, v_1, \ldots, v_r)$ is given by $\sum x^{\alpha_i}e_i <_{\mathbf{w}} \sum x^{\alpha_j}e^j$ if $\sum_i \sum_k \alpha_{i,k} w_k + e_i v_i < \sum_j \sum_k \alpha_{j,k} w_k + e_j v_j$; that is, the order is given by performing a dot product on the combination of the exponents and basis elements with the weight vector.

\begin{lemma}\label{lemmaCancellations}
If $N$ is a submodule of $F$ and $M=F/N$, then the graded Betti numbers of $F/N$ are obtained from those of $F/\mathrm{in}(N)$ by consecutive cancellations. 
\end{lemma}

\proof
By \cite[Proposition 1.8]{Bay82}, there exists a weight vector $(w_0, \ldots w_n, v_1, \ldots, v_m)$ with positive entries such that $x^{\alpha_1}e_i > x^{\alpha_2}e_j$ if and only if $ \alpha_{1,0} w_0 + \cdots + \alpha_{1,n} w_n + f_i v_i > \alpha_{2,0} w_0 + \cdots + \alpha_{2,n} w_n + f_j v_j$. Let $\tilde{N}$ be the homogenization of $N$ as a $\tilde{S} = S[t]$-module with respect to the grading $\deg(x_i) = w_i$ and $\deg(t)=1$ on $\tilde{S}$, and extending to  $\tilde{F} = \oplus \tilde{S}(-f_j)$. Note that every generator of $\tilde{N}$ comes from an element in $N$ with every term other than the initial term multiplied by a positive power of $t$. Also $t$ and $t-1$ are $\tilde{F}/\tilde{N}$-regular. 

Write $\mathbf{F}$ for a minimal graded free resolution of $\tilde{F}/\tilde{N}$. By \cite[Theorem 15.17]{Eis95},  $\mathbf{F} \otimes \tilde{S}/(t)$ is a minimal free graded resolution of $\tilde{F}/\tilde{N} \otimes \tilde{S}/(t) \cong F/\mathrm{in}(N)$, and $\mathbf{F} \otimes \tilde{S}/(t-1)$ is a (not necessarily minimal) graded free resolution of $F/N$. Thus, we can remove a trivial complex from $\mathbf{F} \otimes \tilde{S}/(t-1)$ to get a minimal free resolution of $F/N$, resulting in consecutive cancellation of the Betti numbers from $F/\mathrm{in}(N)$ to $F/N$. $\qed$

\begin{ex}
Consider  $I=(x^2-y^2, xy-z^2)$ and $\mathrm{in}(I)=(x^2, xy, y^3)$ in $S=\CC[x,y,z,w]$, with the graded Betti numbers of $S/I$ 
\[
\beta_{0,0}=1, \ \beta_{1,2} = 2, \ \beta_{2,4}=1
\]
and the graded Betti numbers of $S/\mathrm{in}(I)$
\[
\beta_{0,0}=1, \ \beta_{1,2}=2, \ \beta_{1,3}=1, \ \beta_{2,3}=1, \ \beta_{2,4}=1
\]
differing by the consecutive cancellation of $\beta_{1,3}$ and $\beta_{2,3}$.
\end{ex}

The following construction comes from \cite{Par96}. It will be used in the next lemma.

\begin{defn} \label{def-Pardue}
Let $P = k[z_{ijk}]$ with $0 \leq i \leq n$, $1 \leq k \leq r$ and $1 \leq j < J$
for some $J$ sufficiently large. 
Let $F'$ be the free $P$-module with basis $e'_1, \ldots, e'_r$ with the same generating degrees as $e_1, \ldots, e_r$. 
For a monomial ideal $I$, the \emph{k-polarization} $I^{(p_k)}$ of $I$ is the monomial ideal in $P$ generated by
\[
\{z^{p_k(\mu)} = \prod z_{ijk} | x^{\mu} \text{  is a minimal generator of }I\}.
\]
If $N = I_1 e_1 + \cdots + I_r e_r$ is a monomial submodule of $F$, then the polarization of
$N$ is
\[
N^{(p)} =  I_1^{(p_1)} e'_1 + \cdots + I_r^{(p_r)} e'_r.
\]
Let $L=\{h_{ijk}\}$ be a collection of linear forms in $S$. Define $\sigma_L:P \to S$ by $\sigma_L(z_{ijk}) = h_{ijk}$, and $\sigma'_L:F' \to F$ by $\sigma'_L(\sum f_i e_i) = \sum \sigma_L(f_i) e_i$.
\end{defn}

In particular, given a monomial submodule $N$ and a generic set of linear forms $L$, the operation $\sigma'_L(N^{(p)})$ performs a polarization (that is, adds enough variables so that every term is squarefree), followed by taking generic hyperplane sections. The following lemma shows that the graded Betti numbers are unchanged under $k$-polarization.

\begin{lemma} \label{lem-polarization}
If $N$ is a submodule of $F$ and $L$ is a generic collection of linear forms, then the graded Betti numbers of $F/N$ and $F/\sigma'_L(N^{(p)})$ are the same.
\end{lemma}

\proof See \cite[Corollary 15]{Par96}. $\qed$

\begin{lemma}
If $N$ is a submodule of $F$ and $L$ is the lex submodule of $F$ such that the Hilbert functions of $F/N$ and $F/L$ agree, then $N$ can be deformed to $L$ using a finite sequence of the following operations:
\begin{enumerate}
\item general change of coordinates
\item taking the initial module
\item the operation $\sigma'_L(N^{(p)})$ from Definition \ref{def-Pardue}.
\end{enumerate}
\end{lemma}
\proof See \cite[Proposition 30]{Par96}. $\qed$

\proof
(of Theorem \ref{thmPersistence}):
Given $N$, we deform to the the lex submodule $L$ using the operations above. Clearly the first operation preserves Betti numbers; by Lemma \ref{lem-polarization}, the third operation preserves Betti numbers also. By Lemma \ref{lemmaCancellations}, the second operation can only increase Betti numbers in such a way that the Betti numbers of $F/N$ can be recovered from those of $F/L$ by consecutive cancellations.

By the assumption on $H(F/L, d)$ and $H(F/L, d+1)$, it follows that $L$ has no generators in degree $d+1$. Assume that 
\[
H(M,d+2) < \sum_{i=m-r+1}^m\binom{d+2-f_i+n}{n} + \rho_{d+1}^{\langle d-f_{m-r}+1 \rangle}.
\]
Then $L$ has a generator in degree $d+2$. But $N$ is generated in degrees at most $d$. This must be explained by consecutive cancellation with a non-zero Betti number $\beta_{1,d+2}$ of $L$. 

Since $F/L$ can be written $\oplus_m{S/L_i(-f_i)}$ and lex ideals have linear resolutions by \cite[Proposition 1.3.1]{Cha07}, it follows that the non-zero first syzygy of degree $d+2$ must result from a generator of degree $d+1$, which is a contradiction. $\qed$

Since we have rank-adjusted regularity and persistence, we can use the rank-adjusted Gotzmann number to determine the Grassmannian into which to embed the Quot scheme from \cite{Del15}, $\mathrm{Quot}^P(\mathcal{O}_{\PP^n}^r)$. 

\begin{ex}
Note that the Hilbert polynomial $P(d) = k(d+1) + m$ has Gotzmann number
\[
s=\frac{k(k+1)}{2} + m;
\]
however, when this is being used for $\mathrm{Quot}^P(\mathcal{O}_{\PP^1}^r)$, we can use the rank-adjusted Gotzmann number of $m$. For a specific example, assume that $r=5$ and $k=3$. Then the dimension of the Grassmannian into which the Quot scheme is embedded is as follows: using the standard Gotzmann number, we have a Grassmannian of dimension $(21+4m)(14+4m)$, and using the rank-adjusted Gotzmann number, we have a Grassmannian of dimension $(3+4m)(2+m)$.

The next example shows the improvement in the regularity bound for a rank-$k$ sheaf on $\PP^2$ with Hilbert polynomial 
\[
P(d) = k \binom{d+2}{2} + m_1 (d+1) + m_2.
\]
The standard Gotzmann number is 
\[
s=\frac{1}{24}k(k+1)(3k^2-k+10+12m_1) + \frac{1}{2}m_1(m_1+1) + m_2,
\]
and the rank-adjusted Gotzmann number is
\[
s=\frac{1}{2}m_1(m_1+1) + m_2.
\]
\end{ex}

\section*{Acknowledgements}

The author wishes to thank Vladimir Baranovsky for helpful discussions on this work.

\bibliographystyle{alpha}
\bibliography{Gotzmann_rep}

\end{document}